\documentclass{article}

\usepackage{amsmath}
\usepackage{hyperref}
\usepackage{graphicx}
\usepackage[linesnumbered,ruled,vlined]{algorithm2e}
\usepackage[toc,page]{appendix}
\usepackage{amssymb,amsthm}
\usepackage{authblk}
\usepackage{appendix}
\usepackage{wasysym}
\usepackage{float}
\usepackage{setspace}
\usepackage{enumitem}
\usepackage{bm}
\usepackage{cancel}
\usepackage{pythonhighlight}

\theoremstyle{definition}

\newtheorem{proposition}{Proposition}[subsection]

\theoremstyle{definition}

\setlist[enumerate]{label*=\arabic*.}

\title{Gram quadrature: Numerical integration with Gram polynomials}
\author{Irfan Muhammad}
\affil{School of Computer Science, University of Birmingham, UK}
\begin{document}
	\maketitle

	\begin{abstract}
		The numerical integration of an analytical function $f(x)$ using a finite set of equidistant points can be performed by quadrature formulas like the Newton-Cotes. Unlike Gaussian quadrature formulas however, higher-order Newton-Cotes formulas are not stable, limiting the usable order of such formulas. Existing work showed that by the use of orthogonal polynomials, stable high-order quadrature formulas with equidistant points can be developed. We improve upon such work by making use of (orthogonal) Gram polynomials and deriving an iterative algorithm, together allowing us to reduce the space-complexity of the original algorithm significantly.
	\end{abstract}
	
	
	\section{Introduction} Let us define our problem as follows. We would like to integrate a continuous function $f$ over the interval $[1,1]$, by first approximating $f$ by a polynomial $p$, $$f \approx p$$ and secondly integrating $p$ approximately via a discrete sum (or discretization of the integral), using $N$ samples of $p$ denoted as $p(x_0),p(x_1),\dots,p(x_{N-1})$, each multiplied by an unknown weight parameter $w_k$: $$\int_{-1}^{1}p(x)dx \approx \sum_{i=0}^{N-1}w_ip(x_i)$$
	
	\paragraph{} As for why polynomials can be used for approximations:
	\begin{proposition} For any function $f$ continuous on an interval $[a,b]$, by the Weierstrass Approximation Theorem, there is some polynomial $q$, such that $|f(x)-q(x)| \leq \epsilon$, for all $x \in [a,b]$, where $\epsilon > 0$ and can be made arbitrarily small.
	\end{proposition}\qed
	
	\paragraph{} Assuming $p$ is a polynomial of degree M, its integration can be written as a linear combination of orthogonal polynomials that can characterize  any M dimensional polynomial (or \textit{span} its space). For example, using the orthogonal polynomials $x^0,x^1,x^2,\dots,x^{M-1}$:
	$$    \int_{-1}^1p(x)dx = \sum_{j=0}^{M-1}a_i\int_{-1}^1x^{i}dx$$
	or more generally,
	$$    \int_{-1}^1p(x)dx = \sum_{j=0}^{M-1}a_i\int_{-1}^1\phi_j(x)dx$$
	where $\{\phi_j |j=0,1,\dots,M-1 \}$ is a set of orthogonal polynomials.
	
	\paragraph{} Seeing that we are approximating integrals as a sum, then the above becomes: $$
	\sum_{j=0}^{M-1}a_i\int_{-1}^1\phi_j(x)dx \approx
	\sum_{j=0}^{M-1}a_i \sum_{i=0}^{N-1}w_i\phi_j(x_i)
	$$
	where we have used the same set of weights $\textbf{w}=[w_0,\dots,w_{N-1}]$ and points $\textbf{x}=[x_0,\dots,x_{N-1}]$ for the approximation of each polynomial integration. The above formula gives us a system of linear equations - traditionally known as \textit{normal equations}:
	$$\int_{-1}^1\phi_j(x)dx =\sum_{i=0}^{N-1}w_i\phi_j(x_i)$$
	for all $j=0,1,\dots,M-1$, and $i=0,1,\dots,N-1$. If we have already chosen \textbf{x}, then the system above can be solved for \textbf{w}. That is if let
	$$ b_j = \int_{-1}^1\phi_j(x), \hspace{1cm} A_{i,j} = \phi_j(x_i)$$
	We arrive at the matrix form $\textbf{Aw} = \textbf{b}$, of which is to be solved for \textbf{w}.
	
	\paragraph{} For a certain choice of \textbf{x}, the system above has a unique solution \textbf{w}, and furthermore the approximation is exact. For example, when $N=M$, and the orthogonal polynomials are chosen to be the Legendre polynomials and \textbf{x} is chosen to be the roots of such polynomials, then a solution for \textbf{w} exists and is unique. Furthermore, the integration of each polynomial $\phi_k$ above is done exactly. This is known as \textbf{Gaussian quadrature}.
	
	
	\paragraph{} Note that the quadrature problem is not a \textit{polynomial-fitting} problem, that is
	\begin{flalign*}
		\int_{-1}^{1}f(x)dx &\approx \int_{-1}^{1}p(x)dx &\\
		&=\sum_{i=0}^{M-1}a_i\int_{-1}^1\phi_j(x)dx &\\
		&\approx \sum_{j=0}^{M-1}a_i \sum_{i=0}^{N-1}w_i\phi_j(x_i) &\\
		&=\sum_{i=0}^{N-1}w_ip(x_i) &\\
		&\approx \sum_{i=0}^{N-1}w_if(x_i)
	\end{flalign*}
	and thus not once are the coefficients $\textbf{a}=[a_0,\dots,a_{M-1}]$ approximated or used. Summarily our problem is to compute the integral approximation $$\int_{-1}^{1}f(x)dx \approx \sum_{i=0}^{N-1}w_if(x_i)$$ where $\textbf{w}$ is determined by the normal equations previously.
	
	\paragraph{} In our problem instance, the points \textbf{x} are equidistant within the interval $[-1,1]$, and therefore each $x_i = -1 + 2\frac{i}{N-1}$ for $i=0,1,\dots,N-1$. If the orthogonal polynomials are chosen to be the Lagrange polynomials, we arrive at the \textbf{Newton-Cotes formulas}. However, these formulas are not numerically stable \cite{huybrechs2009stable}, of which showed that numerical stability (to some degree) can be characterized as the positivity of the solution \textbf{w}. That is, as $N\ \&\ M \rightarrow \infty$, the Newton-Cotes solution $\textbf{w}^{\text{NC}}$ for the normal equations is not guaranteed to be strictly positive.
	
	\paragraph{} To derive a stable quadrature rule that is optimal, we improve upon the work of Huybrechs \cite{huybrechs2009stable} who presented a solution via a least-squares approximation to the normal equations. The work however is stated to have been originated from the works of Wilson \cite{wilson1970necessary,wilson1970discrete}. Further developments of Huybrechs' work can be found in \cite{glaubitz2020stable2}, where the restriction on positive weights is removed, and \cite{glaubitz2020stable} where weights are obtained for scattered data points (i.e. not latticed points) on certain multi-dimensional domains. We do not focus on these generalizations however, and leave that for future research.
	
	\section{Least squares quadrature} The least square approximations of the normal equations was presented by Huybrechs in \cite{huybrechs2009stable}. Let us first detail a summary of their solution. 
	
	\begin{proposition} Given an arbitrary linear system $\textbf{Ay} = \textbf{b}$, the least squares approximation of $\textbf{y}$, denoted as $\textbf{y}^*$ (which minimizes  $||\textbf{y}||_2$), is equivalent to solving the system $$\textbf{AA}^{\text{T}}\textbf{u} = \textbf{b}$$ and where $$\textbf{y}^* = A^{\text{T}}\textbf{u}$$ \qed
	\end{proposition}
	
	\paragraph{} Therefore, given the normal equations $\textbf{Aw} = \textbf{b}$, a least-squares approximation of \textbf{w}, denoted as $\textbf{w}^*$, is equivalent to solving the \textit{least-square normal equations}: $\textbf{AA}^{\text{T}}\textbf{u} = \textbf{b}$ for $\textbf{u}$, and then deriving $\textbf{w}^* = \textbf{A}^{\text{T}}\textbf{u} $.
	
	\paragraph{} If we choose a set of polynomials $\bm{\phi}= \{\phi_j |j=0,1,\dots,M-1 \}$, which are pairwise orthogonal with respect to the discrete scalar product over equidistant points \textbf{x}, of which we denote as $u(\phi_k,\phi_l)$, defined as
	$$u(\phi_k,\phi_l) \triangleq \sum_{i=0}^{N-1}\phi_k(x_i)\phi_l(x_i)$$ Then these polynomials have the property that $$u(\phi_k,\phi_l) = 0$$ for all pairs $(\phi_k,\phi_l)$, with each $k,l \in \{0,1,\dots,M-1\}^2$ and $k\neq l$. Using $\bm{\phi}$, we find that in the least-square normal equations, $\textbf{AA}^{\text{T}}$ reduces to a diagonal matrix since $$({AA}^{\text{T}})_{k,l} = \sum_{i=0}^{N-1}\phi_k(x_i)\phi_l(x_i) = 0$$ for all $k,l$ where $k \neq l$. This diagonal matrix can be represented as a single vector $\textbf{v}$ where $v_j \triangleq (\textbf{AA}^{\text{T}})_{j,j}$. With this set of orthogonal polynomials, the equations now become
	$$\textbf{v} \circ \textbf{u} =  \textbf{b}$$
	and therefore
	$$\textbf{w}^* = \textbf{A}^{\text{T}}(\textbf{b}\circ\frac{1}{\textbf{v}} )$$
	\paragraph{}If the polynomials are also orthonormal with respect to the discrete scalar product, then $\textbf{v} \equiv \textbf{1}$. Thus, $\textbf{w}^* = \textbf{A}^{\text{T}}\textbf{b}$.
	
	\paragraph{} The class of polynomials that are orthogonal with the discrete product above, is known in literature as the Gram polynomials\cite{kristinsson1994cross,rivlin1981introduction,barnard1998gram}, or the discrete Legendre orthogonal polynomials \cite{kristinsson1994cross}.
	
	\section{Gram polynomials} We will use the definition for the orthonormal Gram polynomials found in \cite[p. 114]{dahlquist1974numerical}.  Firstly, we will use $N+1$ equidistant points in $[-1,1]$, i.e. $|\textbf{x}| = N+1 $. The maximum polynomial degree is $M$. Let $G_m$ be the $m^{th}$ Gram polynomial, then:
	$$G_{m+1}(x) = \alpha_{m,N}xG_{m}(x) - \gamma_{m,N}G_{m-1}(x)
	$$where
	$$\alpha_{m,N} \triangleq \frac{N}{m+1}\sqrt{\frac{4(m+1)^2-1}{(N+1)^2 - (m+1)^2}} \hspace{1cm} \gamma_{m,N} \triangleq \frac{\alpha_{m}}{\alpha_{m-1}}$$
	and  $$G_0(x) = (N+1)^{-\frac{1}{2}},\hspace{0.2cm}G_{-1}(x)=0 $$
	$$\alpha_{-1} = 1$$
	
	\paragraph{} The required coefficients of these polynomials up to the $M^{th}$ order are $\bm{\alpha} = [\alpha_{-1,N},\alpha_{0,N}, \alpha_{1,N},\dots,\alpha_{M,N}]$ and are \textbf{independent} of the scalar product, unlike the more general and formal theory  \cite{huybrechs2009stable}, which affords a complexity of $O(NM)$ to compute the coefficients. However, the Gram polynomials have a cost of only $O(M)$.
	
	\paragraph{}  Returning to the least-square normal equations, we can now replace each $\phi_i$ with $G_i$ and derive $\textbf{w}^*$.  A unique solution will always exist via a least square approximation \cite{huybrechs2009stable}. Note that in practice, to ensure stability of $\textbf{w}^{*}$, we require that $N \leq M^2$ \cite{huybrechs2009stable,dahlquist1974numerical}.
	
	\subsection{Computing the matrix \textbf{A}} Computing \textbf{A} and storing it  is costly in space. If $M=1000$ and $N = M^2$, then \textbf{A} may take around 8GB of RAM to store. Instead, we learned of a simple idea of computing \textbf{A} only partially from \cite{kristinsson1994cross}. That is, we can compute $\textbf{w}^*$ without storing the entirety of $\textbf{A}$ in memory.
	
	\paragraph{} Firstly, we can rewrite $\textbf{w}^* = \textbf{A}^{\text{T}}\textbf{b}$ as
	\begin{flalign}
		\label{eqn:A_alt_form}
		\textbf{w}^*
		= \sum_{j=0}^{M-1}b_j\textbf{A}_{*,j}^{\text{T}}
		= \sum_{j=0}^{M-1}b_j\textbf{A}_{j,*}
	\end{flalign}
	where we use the notation $\textbf{A}_{*,j}$ to denote the $j^{th}$ column of \textbf{A}, and $\textbf{A}_{k,*}$ the $k^{th}$ row. We will use the rightmost form to compute $\textbf{w}^*$.
	
	\paragraph{} Secondly, realise that due to the definition of the Gram polynomials, the computation of $G_m(x_k)$ for any $x_k$ requires having the computation of only $G_{m-1}(x_k)$ and $G_{m-2}(x_k)$. Therefore, each row of \textbf{A} requires only the previous two rows at most to compute it, i.e. row $j$ only requires row $j-1,j-2$ to compute it. Additionally, seeing that each $A_j$ is used only to be multiplied with $b_j$, we do not need to keep them around after the next two rows $j+1,j+2$ are computed. Hence, combining this fact and eq. \eqref{eqn:A_alt_form}, we can compute $\textbf{w}^*$ by keeping in memory only two rows of \textbf{A}. Therefore the space complexity required reduces to $O(N)$ from $O(NM)$ improving the scalability of the approach significantly. Additionally, realise that nowhere do we actually construct the Gram polynomials. We only require their evaluation at points of \textbf{x}.
	
	\paragraph{} Concretely, let $A_m[n] =  \alpha_{m,N}x_nG_{m}(x_n) - \gamma_{m,N}G_{m-1}(x_n)$. Then,
	
	\begin{flalign*}
		A_{m+1}[n] =  \alpha_{m,N}x_nA_{m}[n] - \gamma_{m,N}A_{m-1}[n]
	\end{flalign*}

	\subsection{Computing the vector \textbf{b}: Moments of Gram polynomials}
	
	\paragraph{} In \cite{huybrechs2009stable}, the author used Gaussian quadrature (GQ) to calculate \textbf{b} exactly, i.e. for each component
	$$b_m = \int_{-1}^{1}G_m(x^{GQ})dx = \sum_{n=0}^{\lfloor\frac{M}{2}\rfloor +1}w^{\text{GQ}}_n G_m(x^{GQ}_n)$$
	where $\textbf{w}^{\text{GQ}}$ and $\textbf{x}^{GQ}$ are GQ weights and points-set using $\lfloor\frac{M}{2}\rfloor +1$ points. Note that for GQ, we can set $N = \lfloor\frac{M}{2}\rfloor +1$, since if $N$ points are used, GQ integrates exactly all polynomials of degree $2N - 1$ or less. From this, one set of GQ points and weights can be used to integrate each of $\{G_m\}_{m=0}^{M}$ exactly. That is, by using GQ with $\lfloor\frac{M}{2}\rfloor +1$ points for the integration of each Gram polynomial.
	
	\paragraph{} Since $G_m$ can be determined recursively, we can compute $\{b_m\}_{m=1}^{M}$ recursively also. Define the vector of points $q_{m}[n] =w^{\text{GQ}}_n G_m(x^{GQ}_n)$ for $n= 1,\cdots,N$. Then,
	\begin{align*}
		q_{m+1}[n] &= w^{\text{GQ}}_n(\alpha_{m,N}x^{GQ}_nG_{m}(x^{GQ}_n) - \gamma_{m,N}G_{m-1}(x^{GQ}_n)) &\\
		&= (\alpha_{m,N} x^{GQ}_nq_{m}[n]) - (\gamma_{m,N} q_{m-1}[n])
	\end{align*}
	and  from this we can obtain $b_{m+1} = \sum_n q_{m+1}[n]$. Initially, we set $q_{-1}[n] = 0$, and $q_{0}[n] = w_n^{\text{GQ}}(m+1)^{-\frac{1}{2}}$. Hence, $b_{-1} = b_{0} =  0$.
	
	\paragraph{} There are algorithms to obtain $\textbf{w}^{\text{GQ}}$ and $\textbf{x}^{GQ}$ (for GQ) with time complexity of $O(M)$ \cite{hale2013fast}. Since each $b_n$ requires $\lfloor\frac{M}{2}\rfloor +1$ summations, the complexity to calculate \textbf{b} is $O(M^2)$.  
	
	\paragraph{} A python implementation for computing $\textbf{w}^*$ is presented in the Appendix.

	\bibliographystyle{abbrv}    
	\bibliography{time_series.bib}
	
	\appendix
	\appendixpage
	\section{Python implementation} We present here a python (3.8) implementation of the Gram quadrature. This requires the python library NumPy (for fast array computations) and Numba (to speed things up).
	\begin{python}                         
		import numpy as np
		from numba import njit

		def create_gram_weights(m):
			"""
			:param m: the degree of the largest possible Gram polynomial for m+1 points.
			:return: gram weights for gram quadrature of m+1 points.
			"""
			
			max_d = int(np.sqrt(m))
			xs = np.linspace(-1, 1, m + 1)
			
			alphas = [1] + [alpha_func(n, m) for n in range(0, max_d + 1)]
			alphas = np.asarray(alphas)
			
			w = np.zeros(m + 1, dtype=float)
			
			# get GQ points and weights
			xs_lg, w_lg = np.polynomial.legendre.leggauss(int(max_d / 2 + 1))
			
			# denote A_p as the previous row of A being worked on. A_n is the newest row.
			A_p = np.zeros(m + 1, dtype=float)
			A_n = np.asarray([(m + 1) ** (-0.5)] * (m + 1))
			
			q_p = np.zeros(int(max_d / 2 + 1), dtype=float)
			q_n = w_lg * np.asarray([(m + 1) ** (-0.5)] * int(max_d / 2 + 1))
			
			# update w for n = 0, 1, 2, 3, ..., max_d
			for j in range(0, max_d + 1):
			w, A_p, A_n, q_n, q_p = update_w(w, A_p, A_n, q_n, q_p, alphas, xs, xs_lg, j)
			
			return w

		@njit
		def alpha_func(n, m):
			a1 = m / (n + 1)
			a2 = (4 * np.power(n + 1, 2)) - 1
			a3 = np.power(m + 1, 2) - np.power(n + 1, 2)
			a4 = a2 / a3
			
			res = a1 * np.sqrt(a4)
		return res

		@njit
		def update_w(w, A_p, A_n, q_n, q_p, alphas, xs, xs_lg, j):
			w += (np.sum(q_n) * A_n)
			
			# update for j = 1, 2, 3, ..., max_d+1
			A_temp = (alphas[j + 1] * (xs * A_n)) - ((alphas[j + 1] / alphas[j]) * A_p)
			A_p = A_n
			A_n = A_temp
			
			q_temp = q_n.copy()
			q_n = (alphas[j + 1] * xs_lg * q_n) - (alphas[j + 1] / alphas[j] * q_p)
			q_p = q_temp
			
			return w, A_p, A_n, q_n, q_p
	\end{python}
	
	As for the usage of the weights then in the same file we can add:
	\begin{python}
		if __name__ == '__main__':
			m = 100
			xs = np.linspace(-1, 1, m + 1)
			
			gram_weights = create_gram_weights(m)
			
			# the sum of stable weights is equal to 2.
			print("Sum of Gram weights:", sum(gram_weights))
			
			# test integration, integrate f below between [-1,1]
			a, b = -1, 1
			f = lambda x: 9 * x ** 2 + 45 * 13 * x ** 3 + 16 * x ** 4
			gram_quad = np.sum(gram_weights * f(xs), axis=-1)
			print("Approx. integration:", gram_quad)
	\end{python}

\end{document}